\documentclass[11pt]{article}

\usepackage{amscd,amsmath, amssymb, fancyhdr, epsfig,url}

\numberwithin{equation}{section}

\newcommand{\version}{version 3.2,\ \   April 1, 2021}

\def\eqref#1{(\ref{#1})}

\newcommand{\arrow}{{\:\longrightarrow\:}}
\newcommand{\Z}{{\Bbb Z}}
\newcommand{\C}{{\Bbb C}}

\newcommand{\R}{{\Bbb R}}
\newcommand{\Q}{{\Bbb Q}}

\def\1{\sqrt{-1}\:}

\newcommand{\cntrct}                % contraction with a vector field
{\hspace{2pt}\raisebox{1pt}{\text{$\lrcorner$}}\hspace{2pt}}

\makeatletter
\def\x@arrow{\DOTSB\Relbar}
\def\xlongequalsignfill@{\arrowfill@\x@arrow\Relbar\x@arrow}
\newcommand{\xlongequal}[2][]{%
        \ext@arrow 0099\xlongequalsignfill@{#1}{#2}}
\def\xlongrightarrowfill@{\arrowfill@\relbar\relbar\longrightarrow}
\newcommand{\xlongrightarrow}[2][]{%
        \ext@arrow 0099\xlongrightarrowfill@{#1}{#2}}
\makeatother

\newcommand{\calo}{{\cal O}}

\renewcommand{\phi}{\varphi}
\renewcommand{\epsilon}{\varepsilon}
\renewcommand{\geq}{\geqslant}
\renewcommand{\leq}{\leqslant}

\newcommand{\diam}{\operatorname{\sf diam}}

\newcommand{\Diff}{\operatorname{Diff}}

\newcommand{\rk}{\operatorname{rk}}

\newcommand{\Kah}{\operatorname{Kah}}

\newcommand{\Teich}{\operatorname{\sf Teich}}
\newcommand{\Proj}{\operatorname{Proj}}
\newcommand{\Pic}{\operatorname{Pic}}
\newcommand{\Comp}{\operatorname{\sf Comp}}

\newcounter{Mycounter}[section]
\newcounter{lemma}[section]
\newcounter{claim}[section]
\renewcommand{\theclaim}{{Claim \thesection.\arabic{claim}}}
\newcommand{\claim}{%
    \setcounter{claim}{\value{Mycounter}}
    \refstepcounter{claim}
    \stepcounter{Mycounter}
    {\noindent \bf \theclaim:\ }}

\newcounter{sublemma}[section]
\newcounter{corollary}[section]
\renewcommand{\thecorollary}{{Corollary \thesection.\arabic{corollary}}}
\newcommand{\corollary}{%
    \setcounter{corollary}{\value{Mycounter}}
    \refstepcounter{corollary}
    \stepcounter{Mycounter}
    {\noindent \bf \thecorollary:\ }}

\newcounter{theorem}[section]
\renewcommand{\thetheorem}{{Theorem \thesection.\arabic{theorem}}}
\newcommand{\theorem}{%
    \setcounter{theorem}{\value{Mycounter}}
    \refstepcounter{theorem}
    \stepcounter{Mycounter}
    {\noindent \bf \thetheorem:\ }}

\newcounter{conjecture}[section]
\newcounter{proposition}[section]
\newcounter{definition}[section]
\renewcommand{\thedefinition}
      {{Definition~\thesection.\arabic{definition}}}
\newcommand{\definition}{%
    \setcounter{definition}{\value{Mycounter}}
    \refstepcounter{definition}
    \stepcounter{Mycounter}
    {\noindent \bf \thedefinition:\ }}

\newcounter{example}[section]
\newcounter{remark}[section]
\renewcommand{\theremark}{{Remark \thesection.\arabic{remark}}}
\newcommand{\remark}{%
    \setcounter{remark}{\value{Mycounter}}
    \refstepcounter{remark}
    \stepcounter{Mycounter}
    {\noindent \bf \theremark:\ }}

\newcounter{problem}[section]
\newcounter{question}[section]
\makeatletter

\@addtoreset{equation}{section} \@addtoreset{footnote}{section} \makeatother

\def\blacksquare{\hbox{\vrule width 5pt height 5pt depth 0pt}}
\def\endproof{\blacksquare}

\begin{document}
%%%%%%%%%%%%%%%%%%%%%%%%%%%%%%%%%%%%%%%%%%%%%%%%%%%%%%%%%%%%
\begin{center}
{\LARGE\bf
Kobayashi pseudometric on hyperk\"ahler manifolds\\[4mm]
}
%%%%%%%%%%%%%%%%%%%%%%%%%%%%%%%%%%%%%%%%%%%%%%%%%%%%%%%%%%%%

 Ljudmila Kamenova, Steven Lu\footnote{Partially supported
 by an NSERC discovery grant},
Misha Verbitsky\footnote{Partially supported by RFBR grants
 12-01-00944-a,  NRI-HSE 
Academic Fund Program in 2013-2014, research grant
12-01-0179, Simons-IUM fellowship, and
AG Laboratory NRI-HSE, RF government grant, ag. 11.G34.31.0023.
 }

\end{center}

%%%%%%%%%%%%%%%%%%%%%%%%%%%%%%%%%%%%%%%%%%%%%%%%
{\small \hspace{0.10\linewidth}
\begin{minipage}[t]{0.85\linewidth}
{\bf Abstract} \\
The Kobayashi pseudometric on a complex manifold is the 
maximal pseudometric such that any holomorphic map from the 
Poincar\'e disk to the manifold is distance-decreasing. 
Kobayashi has conjectured that this pseudometric 
vanishes on Calabi-Yau manifolds. 
Using ergodicity of complex structures, we prove this conjecture
for any hyperk\"ahler manifold  that admits a 
deformation with two Lagrangian fibrations and 
whose Picard rank is not maximal. 
The Strominger-Yau-Zaslow (SYZ) conjecture
claims that parabolic nef line 
bundles on hyperk\"ahler manifolds
are semi-ample. We prove that the Kobayashi 
pseudometric vanishes for any hyperk\"ahler manifold
with $b_2\geq 13$ if the SYZ conjecture holds for all its 
deformations. This proves the Kobayashi conjecture for all K3 surfaces
and their Hilbert schemes. 
\end{minipage}
}
%%%%%%%%%%%%%%%%%%%%%%%%%%%%%%%%%%%%%%%%%%%%%%%%

\tableofcontents

%%%%%%%%%%%%%%%%%%%%%%%%%%%%%%%%%%%%%%%%%%%%%%%%

\section{Introduction}

%%%%%%%%%%%%%%%%%%%%%%%%%%%%%%%%%%%%%%%%%%%%%%%%

The Kobayashi pseudometric on a complex manifold $M$ is the 
maximal pseudometric such that any holomorphic map from the 
Poincar\'e disk to $M$ is distance-decreasing (see 
Section \ref{_kobaya_intro_Subsection_} for 
more details and references). 
Kobayashi conjectured that the Kobayashi pseudometric vanishes 
for all projective varieties with trivial canonical bundle 
(see Problems C.1 and F.3 in \cite{_Kobayashi:1976_}). 
The conjecture was proved for projective K3 surfaces via
the nontrivial theorem in \cite{MM}
that all projective K3 surfaces are swept out by elliptic curves 
(see Lemma 1.51 in \cite{_Voisin:kobayashi_}). We prove the conjecture
for all K3 surfaces as well as for many classes of hyperk\"ahler
manifolds.  %Hyperk\"ahler manifolds naturally generalize K3 surfaces. 
For an extensive survey on problems of Kobayashi and Lang we recommend 
the beautiful survey papers \cite{_Voisin:kobayashi_} by Voisin and 
\cite{_Demailly:kobayashi_} by Demailly. 

Using density arguments and the 
existence of Lagrangian fibrations, it was proved in 
\cite{_Kamenova_V:fibrations_} that all known 
hyperk\"ahler manifolds are Kobayashi non-hyperbolic. Then in 
\cite{_Verbitsky:ergodic_} this result was generalized to all 
hyperk\"ahler manifolds with $b_2 > 3$. All known examples of hyperk\"ahler 
manifolds have  $b_2\geq 7$
and this has been conjectured 
to be true in general. 

We introduce the basics of hyperk\"ahler geometry and
Teichm\"uller spaces in Subsection \ref{_Teich_intro_Subsection_}. 
Upper semicontinuity of the Kobayashi pseudometric is discussed in 
Subsection \ref{_semi-cont_intro_Subsection_}. 
Our main results are in Sections \ref{_hk_Section_} and 
\ref{_Infinitesimal_vanishing_}. 

For a compact complex manifold $M$,
the Teichm\"uller space $\Teich$ is the space of 
complex structures up 
to isotopies. The mapping class group $\Gamma$, or the group of 
``diffeotopies'', acts naturally on $\Teich$. Complex structures with 
dense $\Gamma$-orbits are called {\bf ergodic} 
(see \ref{_ergo_co_Definition_}). We show that if 
the Kobayashi pseudometric on $M$ vanishes, then the Kobayashi 
pseudometric vanishes for all ergodic complex structures on $M$
in the 
same deformation class (\ref{_ergodic_vanishes_Theorem_}). As a corollary, 
the Kobayashi pseudometric vanishes 
for all K3 surfaces (\ref{_K3_vanishing_}), 
and for all ergodic complex structures 
on a hyperk\"ahler manifold  (\ref{_HK_ergodic_vanishing_}). 
%We also prove vanishing of the Kobayashi pseudometric in the 
%non-ergodic case both for those 
%hyperk\"ahler manifolds that admit Lagrangian fibrations. 
%In the latter case, we obtain the stronger result of the
%vanishing of an infinitesimal
%version of the Kobayashi pseudo metric proposed by Royden.
The SYZ conjecture predicts that any hyperk\"ahler manifold has 
a deformation which admits a Lagrangian 
fibration.\footnote{More precisely, SYZ conjecture
is a special form of Kawamata's abundance conjecture
predicting that all parabolic nef line 
bundles on hyperk\"ahler manifolds
are semi-ample; see \cite{_Sawon_} 
and \cite{_Verbitsky:SYZ_} for more details.} 
Assuming this conjecture to be true, we show the vanishing
of the Kobayashi pseudometric for all 
hyperk\"ahler manifolds with $b_2 \geq 13$.
When, in addition, the complex structure is non-ergodic 
%\footnote{This means that its orbit under the diffeomorphism group 
%is dense in the space of all complec structures; see 
%\ref{_ergo_co_Definition_}.}, 
we prove that the infinitesimal Kobayashi pseudometric defined by Royden
vanishes on a Zariski dense open subset of the manifold
(this result is stronger). 

We summarize the main results of this article in the following
theorems; please see the main body of the paper for details
of the definitions and of the proofs.

\hfill

%%%%%%%%%%%%%%%%%%%%%%%%%%%%%%%%%%%%%%%%%%%%%%%%%%%%%%%%%%%%%%
\theorem
Let $M$ be a compact simple hyperk\"ahler manifold 
which satisfies the SYZ conjecture and has a deformation
with non-maximal Picard  rank and two transversal Lagrangian fibrations.
Then the Kobayashi pseudometric on $M$ vanishes.\\

{\bf Proof:} See \ref{_exi_defo_2_fibra_Theorem_} and
\ref{_ergodic_vanishes_Theorem_}. \endproof

\hfill

%%%%%%%%%%%%%%%%%%%%%%%%%%%%%%%%%%%%%%%%%%%%%%%%%%%%%%%%%%%%%%
\theorem
Let $M$ be a compact simple hyperk\"ahler manifold with 
$b_2(M)\geq 13$ which satisfies the SYZ conjecture. 
Then the Kobayashi pseudometric on $M$ vanishes.\\

{\bf Proof:} See \ref{_exi_defo_2_fibra_Theorem_}, 
\ref{_ergodic_vanishes_Theorem_} and
\ref{_two_Lagra_vanishing_Theorem_}. \endproof

\hfill

%%%%%%%%%%%%%%%%%%%%%%%%%%%%%%%%%%%%%%%%
\remark
All known examples of hyperk\"ahler manifolds 
have $b_2(M)\geq 7$ and can be deformed
to one which admits a Lagrangian fibration 
(\cite[Claim 1.20]{_Kamenova_V:fibrations_}).
By the above result, the Kobayashi
pseudometric on known manifolds vanishes,
unless their Picard rank is maximal.

%%%%%%%%%%%%%%%%%%%%%%%%%%%%%%%%%%%%%%%%%%%%%%%%%%%%%%%%%%%
\subsection{Teichm\"uller spaces and hyperk\"ahler
  geometry}
\label{_Teich_intro_Subsection_}
%%%%%%%%%%%%%%%%%%%%%%%%%%%%%%%%%%%%%%%%%%%%%%%%%%%%%%%%%%%

We summarize the definition of the Teichm\"uller space 
of hyper\-k\"ahler manifolds, following \cite{_V:Torelli_}.

\hfill

%%%%%%%%%%%%%%%%%%%%%%%%%%%%%%%%%%%%%%%%%%%%%%%%%%%%%%%%%%%%
\definition
Let $M$ be a compact complex manifold and 
$\Diff_0(M)$ a connected component of its diffeomorphism group
({\bf the group of isotopies}). Denote by $\Comp$
the space of complex structures on $M$, equipped with
a structure of Fr\'echet manifold.
We let $\Teich:=\Comp/\Diff_0(M)$ and call it 
{\bf the Teichm\"uller space} of $M$.

\hfill

\remark
In many important cases, such as in the case
of Calabi-Yau manifolds (\cite{_Catanese:moduli_}), 
$\Teich$ is a finite-dimensional
complex space; usually it is non-Hausdorff.

\hfill

%%%%%%%%%%%%%%%%%%%%%%%%%%%%%%%%%%%%%%%%%%%%%%%%
\definition
Let $\Diff(M)$ be the group of orientable diffeomorphisms of
a complex manifold $M$. Consider the {\bf mapping class group} 
\[ \Gamma:=\Diff(M)/\Diff_0(M)\] acting on $\Teich$. 
The quotient $\Comp/\Diff=\Teich/\Gamma$ 
is called {\bf the moduli space} of complex structures on $M$.
Typically, it is very non-Hausdorff. The set
$\Comp/\Diff$ corresponds bijectively to the set of isomorphism
classes of complex structures.

\hfill

%%%%%%%%%%%%%%%%%%%%%%%%%%%%%%%%%%%%%%%%%%%%%%%%
\definition\label{_hk_Definition_}
A {\bf hyperk\"ahler manifold} is a compact 
holomorphically symplectic manifold admitting a K\"ahler structure.

\hfill

%%%%%%%%%%%%%%%%%%%%%%%%%%%%%%%%%%%%%%%%%%%%%%%%
\definition\label{_hk_simple_Definition_}
A hyperk\"ahler manifold $M$ is called
{\bf simple} if $\pi_1(M)=0$ and $H^{2,0}(M)=\C$.
In the literature, simple hyperk\"ahler manifolds are 
often called {\bf irreducible holomorphic symplectic manifolds},
or simply an {\bf irreducible symplectic varieties}.

\hfill

The equivalence between these two notions is 
based on the following theorem of Bogomolov (via \cite{Yau}) 
that motivated this definition. 

\hfill

%%%%%%%%%%%%%%%%%%%%%%%%%%%%%%%%%%%%%%%%%%%%%%%%
\theorem (\cite{_Bogomolov:decompo_})
Any hyperk\"ahler manifold admits a finite covering
which is a product of a torus and several 
simple hyperk\"ahler manifolds.
\endproof

\hfill

%%%%%%%%%%%%%%%%%%%%%%%%%%%%%%%%%%%%%%%%%%%%%%%%
\remark
Further on, all hyperk\"ahler manifolds
are assumed to be simple, $\Comp$ is the
space of all complex structures of
hyperk\"ahler type on $M$, and $\Teich$
its quotient by $\Diff_0(M)$.

\hfill

%%%%%%%%%%%%%%%%%%%%%%%%%%%%%%%%%%%%%%%%%%%%%%%%
A simple hyperk\"ahler manifold admits a
primitive integral quadratic form 
on its second cohomology group known as 
the Beauville-Bogomolov-Fujiki form. We define it using
the Fujiki identity given in the theorem below;
see \cite{_Fujiki:HK_}. For a more detailed 
description of the form 
we refer the reader to \cite{_Beauville_} and 
\cite{_Bogomolov:defo_}.  

\hfill

%%%%%%%%%%%%%%%%%%%%%%%%%%%%%%%%%%%%%%%%%%%%%%%
\theorem \label{_Fujiki_Theorem_}
(Fujiki, \cite{_Fujiki:HK_}) %\label{_Fujiki_formula_}
Let $M$ be a simple hyperk\"ahler manifold of dimension $2n$ and 
$\alpha\in H^2(M, \Z)$. Then $\int_M \alpha^{2n}=c q(\alpha,\alpha)^n$, 
where $q$ is a primitive integral quadratic form on $H^2(M,\Z)$ of index 
$(3, b_2(M)-3)$, and $c>0$ is a rational number.
\endproof

\hfill

%%%%%%%%%%%%%%%%%%%%%%%%%%%%%%%%%%%%%%%%%%%%%%%%%%%
\remark\label{_Fujiki_gene_Remark_}
Fujiki formula can be used to show that
$\int_M \alpha_1\wedge \alpha_2 \wedge...\wedge \alpha_{2n}$ 
is proportional to a sum of 
$q(\alpha_{i_1}\alpha_{i_2})q(\alpha_{i_3}\alpha_{i_4})...
q(\alpha_{i_{2n-1}}\alpha_{i_{2n}})$ taken over all permutations
$(i_1,i_2,...,i_{2n})$.
Whenever $\alpha,\beta\in H^2(M)$ satisfy $q(\alpha,\alpha)=0$,
Fujiki formula gives 
$$\int_M \alpha^{n}\cup\beta^n=c q(\alpha,\beta)^n.$$

\hfill

%%%%%%%%%%%%%%%%%%%%%%%%%%%%%%%%%%%%%%%%%%%%%%%%%%%%%%
\definition\label{_Fujiki_formula_}
{}From \ref{_Fujiki_Theorem_}, the form $q$ is defined uniquely
up to a sign, except the case of even $n$ and $b_2 \neq 6$. 
To fix the sign, we make the additional
assumption that $q(\omega, \omega)>0$ for every K\"ahler
form $\omega$. Such a form $q$ is called 
{\bf the Bogomolov-Beauville-Fujiki form} (or the 
{\bf BBF form}) of $M$.

\hfill

%%%%%%%%%%%%%%%%%%%%%%%%%%%%%%%%%%%%%%%%%%%%%%%%%
%\remark\label{_q_on_Kahler_remark_}
%The BBF form is remarkably similar to the
%intersection form on second cohomology of a 
%complex surface. In particular,
%for any two K\"ahler classes $\omega, \omega'\in H^2(M,\R)$,
%one has $q(\omega, \omega')>0$ (see e. g. 
%\cite{_Huybrechts:basic_} or \cite{_Boucksom_}).
%
%
%\hfill

The mapping class group
of a hyperk\"ahler manifold can be described
in terms of the BBF form as follows.

\hfill

%%%%%%%%%%%%%%%%%%%%%%%%%%%%%%%%%%%%%%%%%%%%%%%%%%
\theorem
(\cite{_V:Torelli_})
Let $M$ be a simple hyperk\"ahler manifold,
$\Gamma$ its mapping class group, 
and $\Gamma\stackrel \phi \arrow O(H^*(M,\Z), q)$ the natural map. 
Then $\phi$ has finite kernel and its image has finite index
in $O(H^*(M,\Z), q)$.
\endproof

\hfill

%%%%%%%%%%%%%%%%%%%%%%%%%%%%%%%%%%%%%%%%%%%%%%%%%
\definition\label{_monodromy_group_Definition_}
Let $\Teich^I$ be a connected component of the Teichm\"uller
space containing $I\in \Teich$, 
and $\Gamma^I$ the subgroup of the mapping class
group preserving $\Teich^I$. The group  $\Gamma^I$
is called {\bf the monodromy group} of $(M,I)$
(\cite{_Markman:constra_}).

\hfill

%%%%%%%%%%%%%%%%%%%%%%%%%%%%%%%%%%%%%%%%%%%%
\remark\label{_mono_finite_inde_Remark_}
In \cite{_V:Torelli_} it was shown that
$\Gamma^I$ is a finite index subgroup in 
$O(H^*(M,\Z), q)$ independent of $I$.

%%%%%%%%%%%%%%%%%%%%%%%%%%%%%%%%%%%%%%%%%%%%%%%%
\subsection{Ergodic complex structures}
%%%%%%%%%%%%%%%%%%%%%%%%%%%%%%%%%%%%%%%%%%%%%%%%

%%%%%%%%%%%%%%%%%%%%%%%%%%%%%%%%%%%%%%%%%%%%%%%%%%%%%%
\definition\label{_ergo_co_Definition_}
Let $M$ be a complex manifold, $\Teich$ its Teichm\"uller space, and
$I\in \Teich$ a point. Consider the set $Z_I\subset \Teich$
of all $I'\in \Teich$ such that $(M,I)$ is biholomorphic to $(M,I')$. 
Clearly, $Z_I=\Gamma \cdot I$ is the orbit of $I$. A complex structure is 
called {\bf ergodic} if the corresponding orbit $Z_I$ is dense in $\Teich$.

\hfill

%%%%%%%%%%%%%%%%%%%%%%%%%%%%%%%%%%%%%%%%%%%%%%%%%%%%
\theorem
Let $M$ be a simple hyperk\"ahler manifold
or a compact complex torus of dimension $\geq 2$,
and $I$ a complex structure on $M$. Then $I$ is
non-ergodic iff the Neron-Severi lattice of $(M,I)$  
has maximal possible rank. This means that
$\rk NS(M,I)=b_2(M)-2$ for $M$ hyperk\"ahler, and
$\rk NS(M,I)=(\dim_\C M)^2$ for $M$ a torus.

\hfill

{\bf Proof:} See \cite{_Verbitsky:ergodic_}.
\endproof

%%%%%%%%%%%%%%%%%%%%%%%%%%%%%%%%%%%%%%%%%%%%%%%%
\subsection{Kobayashi pseudometric/pseudodistance}
\label{_kobaya_intro_Subsection_}
%%%%%%%%%%%%%%%%%%%%%%%%%%%%%%%%%%%%%%%%%%%%%%%%

Let $M$ be a complex manifold. Recall that a pseudometric on $M$
is a function $d$ on $M\times M$ that satisfies all the properties of a
metric (or distance function) except for the non degeneracy
condition: $d(x,y)=0$ only if $x=y$. The Kobayashi pseudometric 
(a.k.a.~pseudodistance) $d_M$ on $M$ 
is defined as the supremum of all pseudometrics $d$ on $M$ that 
satisfy the distance decreasing property with respect to holomorphic
maps $f$ from the Poincar\'e disk $(\mathbb D,\rho)$ to 
$M$: 
\[ f^*d\leq \rho\ \ \ \mbox{or equivalently}\ \ \ d(f(x),f(y))\leq \rho(x,y)\ 
\, \forall x,y\in {\mathbb D}.
\]
Here $\rho$ denotes the  Poincar\'e metric on ${\Bbb D}$.

The following is S. Kobayashi's standard construction of $d_M$.
Let $$\delta_M(p,q)=\inf \{\, \rho(x,y)\,|\, f: {\mathbb D}\rightarrow M\
\mbox{holomorphic}, f(x)=p, f(y)=q\,\}.$$
Although it does not satisfy the triangle inequality in general,
this is a very useful invariant of the complex structure on 
$M$.  For an ordered subset $S=\{p_1, ..., p_l\}$ of $M$, let
$$\delta^S_M(p,q)=\delta_M(p,p_1)+\delta_M(p_1,p_2)+...
+\delta_M(p_l,q).$$
Then the triangle inequality is attained by setting
$$d_M(p,q)=\inf \delta^S_M(p,q)$$
where the infimum is taken over all finite ordered 
subsets $S$ in $M$.
\\

Royden introduced an infinitesimal version of $d_M$ as follows.
The Ko\-ba\-ya\-shi-Royden Finsler norm on $TM$ is given, for 
$v \in TM$, by 
$$|v|_M=\inf \{\, \frac{1}{R}\,|\, f: {\mathbb D}\rightarrow M\
\mbox{holomorphic}, R>0, f'(0)=Rv\,\}.$$

\noindent
It is the largest ``Finsler" pseudonorm on $TM$ that 
satisfies the distance decreasing property with 
respect to holomorphic maps from the Poincar\'e disk 
and therefore it is automatically ``distance decreasing"
with respect to holomorphic maps.
Royden showed that $|\ |_M$ is upper semicontinuous and
that $d_M$ is the integrated version of $|\ |_M$,
see \cite{Roy}. In particular,
this implies the well known fact that $d_M$ is a continuous 
function for a complex manifold $M$.\\

We recall that both the pseudometric and its 
infinitesimal version are insensitive
to removing complex codimension two subsets of $M$.

\hfill

%%%%%%%%%%%%%%%%%%%%%%%%%%%%%%%%%%%%%%%%%%%%%%%%%%%%%%%%%%%%
\theorem\label{_bira_Kobayashi_Theorem_}
Let $M$ be a complex 
manifold and $Z\subset M$ be a complex analytic 
subvariety of codimension at 
least 2.\footnote{In fact, the same proof would work for
any subset  $Z\subset M$ of Hausdorff codimension at least 3.} 
Then $d_{M\backslash Z}=d_M|_{M\backslash Z}$
and $|\ |_{M\backslash Z}=(|\ |_M)|_{M\backslash Z}$.\\

{\bf Proof:} Theorems 3.2.19 and 3.5.35 in \cite{_Kobayashi:1998_}. 
\endproof

\hfill

%%%%%%%%%%%%%%%%%%%%%%%%%%%%%%%%%%%%%%%%%%%%%%%%%%%%%%%%%%%%
\corollary\label{_bira_Kobayashi_CY_Corollary_}
Let $\tau:\; M\dashrightarrow M'$ be a birational equivalence
of Calabi-Yau manifolds. Suppose that the Kobayashi pseudometric on
$M$ vanishes. Then it vanishes on $M'$.

\hfill

{\bf Proof:} It is easy to check (see subsection 4.4 in 
\cite{_Huybrechts:basic_}) that
the exceptional set of $\tau$ is a subvariety of 
codimension at least 2. Then 
\ref{_bira_Kobayashi_Theorem_} can be applied
to obtain that the Kobayashi pseudometric
vanishes on $M$ and $M'$ (by the distance decreasing
property) whenever it 
vanishes on the smooth locus of $\tau$.
\endproof

%%%%%%%%%%%%%%%%%%%%%%%%%%%%%%%%%%%%%%%%%%%%%%%%
\subsection{Upper semi-continuity}
\label{_semi-cont_intro_Subsection_}
%%%%%%%%%%%%%%%%%%%%%%%%%%%%%%%%%%%%%%%%%%%%%%%%

Recall that a function $F$ on a topological space $X$ with values in 
$\R\cup \{\infty\}$ is upper semi-continuous if and only if
$\{\,x \in X\, |\, F(x) < \alpha\, \}$ is an open set 
for every $\alpha \in \R$. 
It is upper semi-continuous at a point $x_0\in X$ if for all 
$\epsilon>0$ there is a neighbourhood of $x_0$ containing
$\{\,x \in X\, |\, F(x) <  F(x_0)+\epsilon\, \}$. 
If $X$ is a metric space, this is equivalent to 
$$\limsup_{t_i\to t_0} F({t_i})\leq F({t_0}),$$
for all sequence $(t_i)$ converging to $t_0$.
{}From its very definition, 
the infimum of a collection of upper semicontinuous functions
is again upper semicontinuous.\\

We will be interested in the upper semicontinuity of $d_{M_t}$ and
$|\ |_{M_t}$ in the variable $t$ for a proper smooth fibration 
$\pi: {\cal M}\rightarrow T$, i.e., $\pi$ is holomorphic, surjective,
having everywhere of maximal rank and connected fibers
$M_t=\pi^{-1}(t)$. This follows in the standard way as is for the
case of $|\ |_M$ by the following result of Siu.\\

\theorem(\cite{Siu}) 
Let $f:\; D \arrow M$ be a holomorphic immersion of 
a Stein manifold $D$ into a complex manifold $M$. Identify $D$ as
the zero section of the normal bundle $X = f^\ast TM/TD$ 
of $D$ in $M$. 
Then there is a holomorphic immersion of a neighbourhood of
$D$ in $X$ which extends $f$. \endproof\\

\noindent
Since $\pi$ is locally differentiably trivial, we may assume that
$\cal M$ is differentiably a product $T\times M$ and $\pi$ its
projection to the first factor. 
One easily deduce from the above theorem of Siu  applied 
to the graph of a holomorphic map from $D=\mathbb D$ that
$\delta_{J(t)}(p,q)$ and $|v|_{J(t)}$ are upper semicontinuous 
with respect to $p,q\in M$, $v\in TM$ and $t\in T$, where we
have replaced the subscript $M_t$ by its associated complex
structure $J(t)$. It follows then that $\delta^S_{J(t)}(p,q)$ is 
upper semicontinuous with respect to $p,q$ and $t$ and hence
so is $d_{J(t)}(p,q)$. We have established the following
proposition, c.f. \cite{Zai}. \\

\begin{proposition}
Let $\pi: {\cal M}\rightarrow T$ be a proper holomorphic 
and surjective map
having everywhere of maximal rank and connected fibers
$M_t=\pi^{-1}(t)$. Then $d_{M_t}$ and $|\ |_{M_t}$ are upper
semicontinuous with respect to all variables involved, including
$t$. \endproof\\
\end{proposition}

Although we will not need this, a little reflection will show that
one can relax many of the conditions on $\pi$. An immediate
consequence of the above proposition is the following.\\

%%%%%%%%%%%%%%%%%%%%%%%%%%%%%%%%%%%%%%%%%%%%%%%%%%%%%%%
\corollary\label{diam} 
For $M$ a compact complex manifold,
let $\diam(M)$ be the diameter of $M$ with respect to $d_M$.
Then $\diam(M)$ is upper semicontinuous with respect to the
variation of the complex structure on $M$.\\

{\bf Proof:} We need to show that $\diam(M_t)$ 
is upper semicontinuous
with respect to $t$ for a family as given above, i.e. for all $t_0\in T$
and sequences $(t_i)$ converging to $t_0$, 
$$\limsup_{t_i\to t_0} \diam(M_{t_i})\leq \diam(M_{t_0}).$$

If the inequality is false, then after replacing the sequence $(t_i)$ by
a subsequence there is an $\epsilon>0$ such that 
$\diam(M_{t_i})> \diam(M_{t_0})+\epsilon$ for all $i$.
By compactness and the continuity of 
the pseudometric on each $M_t$, 
there exist $p_i, q_i$ such that $\diam(M_{t_i})=d_{M_{t_i}}(p_i,q_i)$.
Replacing by a further subsequence if necessary, we may assume
that the sequences $(p_i)$ and $(q_i)$ are convergent. 
Let $p,q\in M_{t_0}$ be their respective limit. Then by upper
semicontinuity, we have 
$$\diam(M_{t_0})+\epsilon\leq \limsup_{i\to \infty} d_{M_{t_i}}(p_i,q_i)
\leq d_{M_{t_0}}(p,q)\leq \diam(M_{t_0}).$$
This is a contradiction. \endproof

%%%%%%%%%%%%%%%%%%%%%%%%%%%%%%%%%%%%%%%%%%%%%%%%

\section{Vanishing of the Kobayashi pseudometric}% on hyperk\"ahler manifolds}
\label{_hk_Section_}

%%%%%%%%%%%%%%%%%%%%%%%%%%%%%%%%%%%%%%%%%%%%%%%%

%%%%%%%%%%%%%%%%%%%%%%%%%%%%%%%%%%%%%%%%%%%%%%%%%%%%%%%%%%%%
\subsection{Kobayashi pseudometric and ergodicity}
%%%%%%%%%%%%%%%%%%%%%%%%%%%%%%%%%%%%%%%%%%%%%%%%%%%%%%%%%%%%

The main technical result of this paper is the following theorem.
Recall that an ergodic complex structure $I$ on $M$ is one which has
a dense $\Diff(M)$-orbit in the deformation space of complex structures.

\hfill

%%%%%%%%%%%%%%%%%%%%%%%%%%%%%%%%%%%%%%%%%%%%%%%%
\theorem\label{_ergodic_vanishes_Theorem_}
Let $M$ be a complex manifold
with vanishing Kobayashi pseudometric. Then the Kobayashi
pseudometric vanishes for all ergodic complex structures
in the same deformation class.

\hfill

{\bf Proof:} Let $\diam:\; \Teich \arrow \R^{\geq 0}$
map a complex structure $I$ to the diameter of the Kobayashi
pseudodistance on $(M,I)$. By \ref{diam}, this function
is upper semi-continuous. Let $I$ be an ergodic complex structure. 
The set of points $I'\in \Teich$
such that $(M,I')$ is biholomorphic to $(M,I)$
is dense, because $I$ is ergodic. By upper semi-continuity, 
$0=\diam(I) \geq \inf_{I'\in \Teich} \diam (I')$.
\endproof

\hfill

%%%%%%%%%%%%%%%%%%%%%%%%%%%%%%%%%%%%%%%%%%%%%%%%
\corollary\label{_K3_vanishing_}
Let $M$ be a K3 surface. Then the Kobayashi pseudometric
on $M$ vanishes.

\hfill

{\bf Proof:} Notice that any non-ergodic complex
structure on a hyperk\"ahler manifold is projective.
Indeed, if the rank of the Picard group is maximal,
the set of rational $(1,1)$-classes is dense in $H^{1,1}(M)$,
hence the K\"ahler cone contains a rational class and 
$M$ is projective. For all projective $M$, one has $\diam(M)=0$
(see Lemma 1.51 in \cite{_Voisin:kobayashi_} or Corollary 4.5 in
\cite{_Lu:Multiply_}). 
Therefore \ref{_ergodic_vanishes_Theorem_}
implies that $\diam(M)=0$ for non-projective
complex structures as well.
\endproof

\hfill

The same argument leads to the following result.

\hfill

%%%%%%%%%%%%%%%%%%%%%%%%%%%%%%%%%%%%%%%%%%%%%%%%
\theorem \label{_HK_ergodic_vanishing_}
Let $M$ be a hyperk\"ahler manifold admitting
a complex structure with vanishing Kobayashi pseudometric and 
$b_2(M) \geq 4$. 
Then the Kobayashi pseudometric vanishes for
all complex structures $I$ in the same deformation class.

\hfill

{\bf Proof:} The diameter of the 
Kobayashi pseudometric is upper semicontinuous,
by \ref{diam}. Choose any ergodic complex structure
$J$ on $M$ (such $J$ exists because $b_2(M) > 3$). By definition of 
ergodic complex structures, 
in any neighbourhood of $I$ one has
a complex manifold isomorphic to $(M,J)$. By upper semicontinuity, 
one has $\diam (M,J)\leq \diam (M,I)=0$. Now vanishing 
of the Kobayashi pseudometric follows from 
\ref{_ergodic_vanishes_Theorem_}.
\endproof

%%%%%%%%%%%%%%%%%%%%%%%%%%%%%%%%%%%%%%%%%%%%%%%%%%%%%%%%%%%%%%%%%%%%
\subsection{Lagrangian fibrations in hyperk\"ahler geometry}
%%%%%%%%%%%%%%%%%%%%%%%%%%%%%%%%%%%%%%%%%%%%%%%%%%%%%%%%%%%%%%%%%%%%

The theory of Lagrangian fibrations on hyperk\"ahler manifolds
is based on the following remarkable theorem of D. Matsushita 

\hfill

%%%%%%%%%%%%%%%%%%%%%%%%%%%%%%%%%%%%%%
\theorem\label{_Matsushita_Theorem_}
(\cite{_Matsushita:fibred_})
Let $M$ be a simple hyperk\"ahler manifold, and
$\phi:\; M \arrow X$ a surjective holomorphic map, with
$0<\dim X < \dim M$. 
Then the fibers of $\phi$
are Lagrangian subvarieties on $M$, and the general
fibers of $\phi$ are complex 
tori.\footnote{These fibers are known to be abelian varieties,
%even if torus is non-algebraic; 
see \cite[Proposition 3.3]{_COP:non-alge_}.}
\endproof

\hfill

%%%%%%%%%%%%%%%%%%%%%%%%%%%%%%%%%%%%%%%%%%%%%%%%%%%%%%%%%%%
\remark
Such a map is called
{\bf a Lagrangian fibration}. 
All the known examples of hyperk\"ahler
manifolds admit Lagrangian 
fibrations (\cite[Claim 1.20]{_Kamenova_V:fibrations_}).

\hfill

%%%%%%%%%%%%%%%%%%%%%%%%%%%%%%%%%%%%%%%%%%%%%%
\definition
 A cohomology class  $\eta\in H^2(M,\R)$ is called
{\bf nef} if it lies in the closure of the K\"ahler cone;
a line bundle $L$ is nef if $c_1(L)$ is nef.
A nef line 
bundle $L$ is {\bf big} if $\int_M c_1(L)^{\dim_\C M}\neq 0$. 
A non-trivial nef line bundle $L$ on a hyperk\"ahler manifold
is called {\bf parabolic} if it is not big. 
{}From the definition of the BBF form, this is 
equivalent to $q(c_1(L), c_1(L))=0$.
Lagrangian fibrations are in 
bijective correspondence with semiample parabolic bundles, 
as follows from Matsushita's theorem. 

\hfill

%%%%%%%%%%%%%%%%%%%%%%%%%%%%%%%%%%%%%%%%%%%%%%%%%%%%%%%%%%%
\claim\label{_Lagra_nef_Claim_}
Let $M$ be a simple hyperk\"ahler manifold, and
$L$ a non-trivial semiample bundle on $M$. Assume that
$L$ is not ample. Consider the holomorphic map
$\pi:\; M \arrow \Proj(\bigoplus_N H^0(M,L^N)$.
Then $\pi$ is a Lagrangian fibration.
Moreover, every Lagrangian fibration is uniquely determined by a 
parabolic nef line bundle.

\hfill

{\bf Proof:} The first statement of \ref{_Lagra_nef_Claim_} is a corollary 
of \ref{_Matsushita_Theorem_}. 
Let  $M \stackrel \pi \arrow X$ be a Lagrangian fibration. By Matsushita's 
results  (\cite{_Matsushita:fibred_}), $X$ is projective and 
$H^*(X)\cong H^*(\C P^n)$. 
Denote by $\eta \in H^2(M,\Z)$ the ample generator
of $\Pic(X)$. Then $\pi^*\eta=c_1(L)$, where $L=\pi^*\calo_X(1)$
is a parabolic nef bundle on $M$.
\endproof

\hfill

The {\bf SYZ conjecture} (\cite{_Sawon_}, \cite{_Verbitsky:SYZ_})
claims that any parabolic nef line bundle 
on a hyperk\"ahler manifold is semiample, that is, it is
associated with a Lagrangian fibration. This is true 
for K3 surfaces (as it follows from the Riemann-Roch formula) 
and for all deformations of Hilbert schemes of K3 surfaces 
(\cite{_Markman:SYZ_} and \cite{_Bayer_Macri_}). 

\hfill

Further on, we shall need the following birational version of 
Matsushita's theorem on Lagrangian fibrations, which is 
due to Matsushita-Zhang.

\hfill

%%%%%%%%%%%%%%%%%%%%%%%%%%%%%%%%%%%%%%%%%%%%%%%%
\theorem\label{_Matsushita_Zhang_nef_Theorem_}
(\cite[Theorem 1.4]{_Matsushita_Zhang_})
Let $X$ be a projective hyperk\"ahler manifold, and
$\overline {BK}(X)$ be the closure of a 
union of all K\"ahler cones for all birational models of $X$.
Conisder an effective $\R$-divisor $P\in \overline {BK}(X)$. Then 
there exists a birational modification
$\tau:\; X' \dashrightarrow X$, where $X'$ is 
a projective hyperk\"ahler manifold such that $\tau^* P$ is nef.
\endproof

\hfill

%%%%%%%%%%%%%%%%%%%%%%%%%%%%%%%%%%%%%%%%%%%%%%%%
\theorem\label{_bira_lag_Theorem_}
Let $M$ be a projective hyperk\"ahler manifold, and
$L$ a line bundle of Kodaira dimension $\frac 1 2 \dim_\C M$.
Then there exists a birational modification
$\tau:\; M' \dashrightarrow M$ from a projective hyperk\"ahler
manifold such that $\tau^* L$ is semiample, and induces
a Lagrangian fibration as in \ref{_Lagra_nef_Claim_}.

\hfill

{\bf Proof:} 
Let $L$ be a nef bundle on a K\"ahler manifold.
Recall that the {\bf numerical Kodaira 
dimension} of $L$ is the maximal
$k$ such that $c_1(L)^k\neq 0$. The {\bf Kodaira dimension}
of $L$ is the Krull dimension of the ring $\bigoplus_N H^0(M,L^N)$.

Consider the modification
$\tau:\; M' \dashrightarrow M$ produced by the
Matsushita-Zhang theorem. Then the numerical dimension
of $\tau^* L$ is equal to $\frac 1 2 \dim_\C M$,
by \cite{_Verbitsky:coho_announce_}, and the
Kodaira dimension stays the same. As shown in
\cite[Theorem 1]{_Kawamata:pluricanonical_} (see also
\cite[Proposition 2.8]{_BCEKPRSW_}),
whenever the numerical dimension of a nef bundle 
is equal to its Kodaira dimension, the bundle is semiample.
Then \ref{_bira_lag_Theorem_} follows from 
\ref{_Lagra_nef_Claim_}. \endproof

\hfill

This result motivates the following definition.

\hfill

%%%%%%%%%%%%%%%%%%%%%%%%%%%%%%%%%%%%%%%%%%%%%%%%%%
\definition
Let $\tau:\; M' \dashrightarrow M$  
be a birational map of hyperk\"ahler manifolds,
and ${\cal L}$ a Lagrangian fibration on $M$.
Then $\tau^* {\cal L}$ is called {\bf a birational
Lagrangian fibration on $M'$}. Its fibers are
proper preimages of those fibers of ${\cal L}$ which
are not contained in the exceptional locus of $\tau$.

%%%%%%%%%%%%%%%%%%%%%%%%%%%%%%%%%%%%%%%%%%%%%%%%%%%%%%%%%%%%%%%%%%%%
\subsection{Kobayashi pseudometrics and Lagrangian fibrations}
%%%%%%%%%%%%%%%%%%%%%%%%%%%%%%%%%%%%%%%%%%%%%%%%%%%%%%%%%%%%%%%%%%%%

The idea to use \ref{_two_Lagra_vanishing_Theorem_}
is suggested by Claire Voisin. We are very grateful 
to Prof. Voisin for her invaluable help.

\hfill

%%%%%%%%%%%%%%%%%%%%%%%%%%%%%%%%%%%%%%%%%%%%%%%%%%%%%%%
\theorem\label{_two_Lagra_vanishing_Theorem_}
Let $M$ be a simple hyperk\"ahler manifold admitting two 
Lagrangian fibrations associated with two non-proportional
parabolic classes. Then the Kobayashi pseudometric on $M$ vanishes.

\hfill

{\bf Proof:} Let $\pi_i:\; M \arrow X_i, \; i=1,2,$
be the Lagrangian fibration maps. Since the general
fibers of $\pi_i$ are tori, the Kobayashi pseudometric
vanishes on each fiber of $\pi_i$. To prove that 
the Kobayashi pseudometric vanishes on $M$, it would suffice
to show that a general fiber of $\pi_1$ intersects
all the fibers of $\pi_2$. 

Let now $\omega_i$ be an ample class of $X_i$ lifted
to $M$, and $2n=\dim_\C M$. Since $\omega_1$ and $\omega_2$  are not
proportional, the standard linear-algebra argument,
often called the Hodge index formula, implies that 
$q(\omega_1, \omega_2) \neq 0$. Indeed, 
$q(\omega_1, \omega_2)\neq 0$ or else
the space $(H^{1,1}(M,\R), q)$ would contain a 
2-dimensional isotropic plane while
its signature is $(1, b_2-3)$. 

Clearly, the fundamental class $[Z_i]$ 
of a fiber of $\pi_i$ is proportional to $\omega_i^n$. Fix the
constant multiplier in such a way that $[Z_i]=\omega_i^n$. 
The fibers of $\pi_1$ intersect 
that of $\pi_2$ if $\int_M [Z_1]\wedge [Z_2]>0$.
However, Fujiki's formula (see \ref{_Fujiki_gene_Remark_}) shows that 
$\int_M [Z_1]\wedge [Z_2]=Cq(\omega_1,\omega_2)^n>0$.
This means that $Z_1$ and $Z_2$ always intersect.
\endproof

\hfill

To prove that a given hyperk\"ahler manifold admits
a deformation with two distinct Lagrangian fibrations,
we use an argument based on \cite{_AV:construction_autom_}.

\hfill

%%%%%%%%%%%%%%%%%%%%%%%%%%%%%%%%%%%%%%%%%%%%%%%%%%%%%%%%%%%%
\theorem\label{_manifold_with_parabo_AV_Theorem_}
Let $M$ be a simple hyperk\"ahler manifold
with $b_2(M) > 13$. Then $M$ admits a projective deformation
with Picard rank equal to 5, with the K\"ahler cone
equal to its positive cone, and its automorphism group
has finite index in the arithmetic group $SO(\Pic(M))$
of orthogonal automorphisms of its Picard lattice.

\hfill

{\bf Proof:} From \cite[Theorem 3.11]{_AV:construction_autom_}
it follows that there exists a projective deformation with
Picard rank 5 and without MBM classes of type (1,1). 
From \cite[Theorem 2.10]{_AV:construction_autom_}
it follows that for such a manifold the  K\"ahler cone
equal to its positive cone, and from 
\cite[Theorem 2.6, Theorem 2.7, Corollary 2.12]{_AV:construction_autom_}
it follows that its automorphism group
has finite index in the arithmetic group $SO(\Pic(M))$.
\endproof

\hfill

%%%%%%%%%%%%%%%%%%%%%%%%%%%%%%%%%%%%%%%%%%%%%%%%%%%%%%%%%%%%
\theorem\label{_exi_defo_2_fibra_Theorem_}
Let $M$ be a projective, simple hyperk\"ahler manifold
with Picard rank $\geq 5$ and with the K\"ahler cone
equal to its positive cone. Assume that $M$ satisfies the
SYZ conjecture, that is, any nef bundle on
$M$ is semiample. Then $M$ admits infinitely many
transversal holomorphic Lagrangian fibrations.
In particular, the Kobayashi pseudometric on $M$ vanishes. 

\hfill

{\bf Proof:} By Meyer's Theorem (\cite{_Meyers_}),
any indefinite quadratic lattice of indefinite
signature and $\rk \geq 5$ represents 0. Therefore,
there exist rational vectors on the boundary of the
positive cone (hence, the K\"ahler cone $\Kah(M)$) of $M$.
In fact, since the automorphism group of 
$M$ is arithmetic, it acts on the set of such vectors
with infinite orbits, hence the set of the rational
vector on the boundary of $\Kah(M)$ is infinite. 
Each such vector corresponds to a Lagrangian fibration,
because we assume that the SYZ conjecture holds. 
\endproof

\hfill

\remark
The SYZ conjecture is true for all known hyperk\"ahler examples, i.e., 
for deformations of Hilbert schemes of points on K3 surfaces 
(Bayer-Macr\`i \cite{_Bayer_Macri_}; Markman \cite{_Markman:SYZ_}), 
for deformations of the generalized Kummer varieties (Yoshioka 
\cite{_Yoshioka_}), for O'Grady's sixfolds (Mongardi-Rapagnetta 
\cite{_Mongardi_Rapagnetta_}), and for O'Grady's tenfolds (Mongardi-Onorati, 
\cite{_Mongardi_Onorati_}).

\section{Vanishing of the infinitesimal pseudometric}
\label{_Infinitesimal_vanishing_}

In this section, 
we are interested in conditions that guarantee the vanishing 
of the infinitesimal Kobayashi pseudometric 
$|\ |_M$ on a Zariski dense open subset of $M$. Recall 
that the SYZ conjecture predicts the existence of a Lagrangian
fibration for every hyperk\"ahler 
$M,\ \dim_\C M=2n$. Furthermore, if 
the base of the fibration 
is smooth (this is conjectured, see \cite{Hwang}),
then the base is isomorphic to $\C P^n$, as shown by Hwang
(see \cite{Hwang, GL}).
If $M$ is projective and admits an
abelian fibration, then we have the following 
two results.\\

\theorem\label{abelian_f} Let $M$ be a projective manifold with an 
equidimensional abelian fibration
$f:M\to B$ (holomorphic surjective 
with all fibres of the same dimension
and general fibres isomorphic to abelian varieties) 
where $B$ is a complex projective
space of lower dimension. If $f$ has no multiple fibres 
in codimension one, then $|\ |_M$ vanishes 
everywhere on $M$. In particular, if $M$ is a projective 
hyperk\"ahler manifold with a birational Lagrangian fibration
over a nonsingular base
without multiple fibres in codimension one, 
then $|\ |_M$ vanishes everywhere.\\

{\bf Proof}: Let $v\in T_xM$. Then $v$ can be regarded as
the first order part of some non vertical $k$-jet $\nu$, which
in turn push forward to a non-trivial jet prescription $\mu$ at
$b=f(x)\in B$. This jet prescription $\mu$ is clearly satisfied by
an algebraic holomorphic map $h:\C\to B$. 
Since this map can be chosen to avoid any subset of 
codimension two or more, we see by so doing that the pull back 
fibration $M_h\to \C$ has no multiple fibres. Hence  
all higher order jet infinitesimal pseudometric vanishes
on $M_h$ by \ref{jet-vanishing}. Since $\nu$ is in the 
image of a $k$-jet on $M_h$, it also has zero $k$-jet infinitesimal
pseudometric by the distance decreasing property and therefore 
$|v |_M=0$. \endproof\\

\theorem\label{_abelian_} Let $M$ be a projective manifold. 
Let $f:M\rightarrow B$ 
define an abelian fibration. Assume that 
%$B$ is birational to projective space ${\mathbb P}^k$ or that 
there is a subvariety
$Z\subset X$ that dominates $B$ and is birational to
an abelian variety. Then
$|\ |_M$ vanishes everywhere above a Zariski dense open subset 
$U$ in $B$. In particular, this holds for hyperk\"ahler 
manifolds with  
$b_2 \geq 5$ having two birational Lagrangian fibrations.\\

{\bf Proof:} 
By hypothesis and the resolution of singularity theorem,
$Z$ is the holomorphic image of a nonsingular projective 
variety $A$ obtained from an abelian variety by blowing
up smooth centres. By construction, any vector in $A$ is 
in the tangent space of an entire holomorphic curve. Let
$g:A\to B$ be the composition with the projection to $B$
and $\mbox{disc}(g)$ its discriminant locus. Let $v\in TM$
be a nonzero vector above the complement
$U$ in $B$ of $\mbox{disc}(f)\cup\mbox{disc}(g)$. 
If $v$ is vertical, then
it is a vector on the fibre $A$ through $p$, which
is an abelian variety and clearly $|v|_M\leq |v|_A=0$
in this case. If $v$ is horizontal, then there is a vector $v'$
in $TA$ by construction such that $f_*v=g_*v'$. Let 
$h:\C\to A$ be such that $h'(0)=v'$ and $\pi:M_h\to \C$ 
be the pull back fibration via the base change by $g\circ h$.
Then $\pi$ has no multiple fibres and $v$ lies in the
image of $TM_h$ by construction. \ref{vanishing} from Appendix
now applies to show that $|\ |_{M_h}$ vanishes and so
the distance decreasing property yields $|v|_M=0$.
The last statement follows from the projectivity of  $M$
by the proof of \ref{_K3_vanishing_}.
\endproof\\

%\corollary
%Let $M$ be a projective hyper\"ahler manifold.
%Assume that $M$ has a 
%Lagragian fibration (i.e., the SYZ conjecture 
%holds for $M$) and that the base of this fibration 
%is a projective space. Then $|\ |_M$ vanishes 
%everywhere above a dense Zariski open subset
%of the base. \endproof\\
%
%It is expected that Lagrangian fibrations on 
%hyperk\"ahler manifolds are holomorphic
%and with projective spaces as bases and
%it is known that the base of a holomorphic Lagrangian 
%fibration from a hyperk\"ahler manifold is a projective
%space if the base is smooth, see \cite{GL}.
%In order to apply this theorem to Lagragian fibrations,

%%%%%%%%%%%%%%%%%%%%%%%%%%%%%%%%%%%%%%%%%%%%%%%%%%%%%%%%%%%%

\section{Appendix on abelian fibrations}

The following are some relevant basic results
concerning abelian fibrations found in 
\cite{_Lu:Multiply_}, which was cited and used in \cite{_Camp04_},
\cite{_Lu:On_Kob_} and \cite{_Lu:Hol_}. Recall that a fibration is
a proper surjective map with connected fibres. 
All fibrations are assumed to be projective in this section
and abelian fibrations are those whose general fibres are
abelian varieties.\\

\begin{proposition}\label{BLT}
Let $e:P \rightarrow \mathbb{D}$ define an abelian fibration which, outside $0\in \mathbb{D}$, 
is smooth with abelian varieties  as fibers. Let $n_0$ be the 
multiplicity of the central fiber $P_0$. Then there is a component of 
multiplicity $n_0$ in $P_0$.\\
\end{proposition}

%%\noindent
{\bf Proof} We may reduce the problem to  the case of $n_0=1$ by
the usual base change $z\mapsto z^{n_0}$ so that the resulting 
object (after normalization) is again such a fibration with an 
unramified cover to the original $P$. Let $\{m_1,m_2, ..., m_k\}$
be the set of multiplicities of the 
components of $P_0$. By assumption,
there exists integers $l_i$ such that $l_1m_1+l_2m_2+...+l_km_k=1$.

As fibrations are assumed to be projective in this paper, we may
assume that $f$ is algebraic.
By restricting to $\mathbb{D}_r=\{\, t\, :\ |t|<r\,\}$ for an $r<1$ if necessary, 
we can construct an algebraic multi-section
$s_i$ with multiplicity $m_i$ above $\mathbb{D}_r$ by simply taking a
one dimensional algebraic slice transversal to the $i$-th component for each
$i$.  Above each point $t$ outside $0$, $s_i$ consists of $m_i$ points
$s_i^j(t),\ j=1,2,...,m_i$. Then it is easy to verify that 
$$(l_1\sum_j s^j_1(t)) + (l_2\sum_j s^j_2(t)) + ... + (l_k\sum_j s^j_k(t))$$
is independent of the choice of an origin in the abelian variety $P_t$.
This gives a section $s$ of the fibration outside $0$ and we now show that
$s$ must be algebraic, giving a section of $f$ and establishing our 
proposition. 

	This can be accomplished by looking at the base change
via $z\mapsto t=z^m$ where $m$ is the least common multiple of 
$m_1,m_2, ..., m_k$. Then each $s_i$ lifts to $m_i$ sections which the 
cyclic Galois action permutes. Hence, the Galois action
of $Z_m$ acts transitively on the $m$ sections constructed by replacing 
the $i$th term in the above expressing with each of the $m_i$ sections
and so this set of sections descends to a section of the original fibration
as desired. 
\endproof \\

%\noindent
We remark that 
the above proposition is really a special case of a result of Lang and
Tate found in \cite{LT}. \\

This proposition allows us to do exactly the 
same analysis as in the case for elliptic fibrations 
done in \cite{BL} to
obtain the following theorem. 
We refer the reader 
there or to \cite{_Lu:Multiply_} for the 
detail of the proof.\\

\theorem\label{vanishing}
Let $f:X\rightarrow C$ 
define an abelian fibration over a complex curve $C$. Then,
for each $s\in C$, the multiplicity of the fibre $X_s$ at $s$ 
is the same as the minimum multiplicity $m_s$ of the components
of $X_s$. 
Let  the $\Q$-divisor $A=\sum_s(1-1/m_s)s$ 
be the  resulting orbifold structure on $C$. 
 Then the three conditions $d_X= 0$ on $X$, 
$|\ |_{X}= 0$ on $X$
and $(C,A)$ is nonhyperbolic $($that is, $C$ is quasiprojective and
$e(C)-\deg A\geq 0)$ are equivalent for such a fibration. In the case
$C$ is quasiprojective, these three conditions are equivalent to
the absence of non-commutative free subgroups in $\pi_1(X)$ and to
$\pi_1(X)$ being solvable up to a finite extension.\\

{\bf Proof:} In the case $(C,A)$ is uniformizable, we 
may pull back the fibration to the universal cover $U$ of $(C,A)$
with resulting fibration $\tilde f: Y\to X$. This is the case
when $C$ is not quasi projective and otherwise when 
$e(C)\leq \deg A$, with equality if and only if $U=\C$, and 
when $e(C)>\deg A$ in which case either $C=U=\C$ and 
$A$ is supported at one point, or $C={\mathbb P}^1$ and
$A$ is supported at more than two points, see for example
\cite{FK}. 
In all these cases
$U$ is non hyperbolic if and only if $(C,A)$ is.
By construction $Y$ has no multiple
fibres over $U$ and is unramified over $X$
(in codimension one) so that all holomorphic curves in $X$
lifts to $Y$. Hence the Kobayashi pseudometrics and
norms vanish on $X$ if and only if it is so on $Y$ and so we only
need to show the vanishing of $|\ |_Y$ in this case since the
fundamental group characterization in the quasi projective 
case follows from the same characterization of the Galois
group of the uniformization $\tilde f:U\to C$ and the exact
sequence of fundamental groups of a fibration without 
multiple fibres. Note that in the case $U={\mathbb P}^1$,
to show that $|\ |_Y$ vanish at a point above $z\in U$
we may replace $U$ by $\C=U\setminus \{z\}$ since
$|\ |_Y\leq |\ |_{Y'}$ where 
$Y':=Y\setminus \tilde f^{-1}(z)\subset Y$.\\
\indent In the case
$(C,A)$ is not uniformizable, then $C={\mathbb P}^1$
and $A$ is supported at one or two points and the exact
sequence of orbifold fundamental groups shows that
$\pi_1(X)$ is a quotient of $\pi_1(X_s)$ for a general
fibre $X_s$, hence abelian. Thus it suffice to show that
for a point $p\in X$, $|\ |_X$ vanishes there in this case
and for this it is sufficient to replace $X$ by the complement
of a fibre $X_z$ different from the fibre 
$X_w$ containing $p$ and $C$ by $C\setminus \{z\}$, 
where in the case $w$ lies in the 
support of $A$, we choose $z$ to be the other point in this
support if one exists. Then $(C,A)$ is uniformizable by $\C$.\\
\indent Hence it remains to show that $|\ |_X$ vanishes at a
point $p$ for the case $X$ has no multiple fibres and $C=\C$.
In fact, given a finite jet prescription at $p$, we can find an entire
holomorphic curve through $p$ satisfying the jet prescription as
follows. The jet prescription gives rise to a jet prescription at 
$f(p)\in C$ which we assume without loss of generality to 
be the origin of $C=\C$. Let $l$ be the first non vanishing order
of the latter jet and let $\tilde f:Y\to \C$ be the pull back fibration 
by the base change $z\mapsto z^l$. Then the inverse function
theorem allows us to translate the jet prescription at $p$ to a
section jet prescription on $Y$ over $0\in \C$. As there are 
no multiple fibres for $\tilde f$, 
proposition~\ref{BLT} yield the existence of local
sections of $\tilde f$ through any point of $\C$. 
The Cousin principles apply in this
situation (i.e., an analogue of Weierstrass' theorem can be
worked out, see \cite{BL,_Lu:Multiply_}) so that 
we can patch up a minimal covering
family of such sections, including the one with the jet prescription,
to give a global section of $\tilde f$ with the jet prescription and
this gives the required entire holomorphic curve. \endproof\\

Instead of restricting our attention 
to just the first order jets for the infinitesimal 
pseudometric, one can generalize the definition of $|\ |_X$ to
jets of arbitrary finite order, see \cite{_Lu:Multiply_}.
By their very definition, these infinitesimal 
pseudometrics dominates
$|\ |_X$ by truncating the jets 
to first order. The exact same proof as
above yields the following generalization, see \cite{_Lu:Multiply_}.\\

\theorem\label{jet-vanishing} The above theorem holds if
$|\ |_X$ is replaced by its more general $k$-th order jet
version, for all integer $k>0$.\\

\hfill

{\bf Acknowledgements:} We would like to express our special gratitude to 
Claire Voisin. This project started when the authors were visiting CRM in 
June 2013. The first named author would like to thank BICMR, where she worked 
on this paper during her visit. The authors also thank Eyal Markman, 
Rob Lazarsfeld, Igor Dolgachev, Stephane Druel and Kieran O'Grady for their 
answers and comments. The first and the last named authors thank the SCGP 
where they enjoyed working and completing this paper. We would like to 
thank the referee for the useful suggestions. Many thanks to Christian Lehn
who told us about an error in the published version of this paper
(the proofs of Theorem 2.12, and Theorem 2.13 were wrong) 
which are now replaced by corrected arguments. 

%\hfill

\newpage

{\small

\noindent {\sc Ljudmila Kamenova\\
Department of Mathematics, 3-115 \\
Stony Brook University \\
Stony Brook, NY 11794-3651, USA,} \\
\tt kamenova@math.sunysb.edu
\\
 
\noindent {\sc Steven Lu\\
D\'epartment de Math\'ematiques, PK-5151\\
Universit\'e du Qu\'ebec \`a Montr\'eal (UQAM)\\
C.P. 8888 Succersale Centreville H3C 3P8,\\ 
Montr\'eal, Qu\'ebec, Canada,} \\
\tt lu.steven@uqam.ca\\

\noindent {\sc Misha Verbitsky\\
{\sc Laboratory of Algebraic Geometry,\\
National Research University HSE,\\
Faculty of Mathematics, 7 Vavilova Str. Moscow, Russia,}\\
\tt  verbit@mccme.ru}, also: \\
{\sc Kavli IPMU (WPI), the University of Tokyo.}
}


\begin{thebibliography}{AV1}

%\bibitem[AC]{_Amerik_Campana:2008_}
%E. Amerik, F. Campana,
%{\em Fibrations m\'eromorphes sur certaines vari\'et\'es 
%\`a fibr\'e canonique trivial,} 
%Pure Appl. Math. Q. 4 (2008), no. 2, part 1, 509-545.

\bibitem[AV]{_AV:construction_autom_} 
E. Amerik, M. Verbitsky,
{\em   Construction of automorphisms of hyperk\"ahler
  manifolds}, Compos. Math. 153 (2017), no. 8, 1610-1621.


\bibitem[BCEKPRSW]{_BCEKPRSW_} 
T. Bauer, F.Campana, T. Eckl, S. Kebekus, T. Peternell, S. Rams, 
T. Szemberg and L. Wotzlaw, 
{\em A reduction map for nef line bundles,} 
Complex geometry (G\"ottingen, 2000) Springer, 
Berlin, 27-36.

\bibitem[BM]{_Bayer_Macri_}
Bayer, A., Macri, E.
{\em MMP for moduli of sheaves on K3s via wall-crossing: nef and movable 
cones, Lagrangian fibrations,} arXiv:1301.6968 [math.AG].

\bibitem[Bea]{_Beauville_} 
 Beauville, A. {\em 
Varietes K\"ahleriennes dont la premi\`ere classe de Chern est
nulle.}  J. Diff. Geom. {\bf 18} (1983) 755-782.

%\bibitem[Bes]{_Besse:Einst_Manifo_} 
%Besse, A., {\em Einstein Manifolds}, Springer-Verlag, New York (1987)

\bibitem[Bo1]{_Bogomolov:decompo_}  
Bogomolov, F. A., {\em On the decomposition of 
K\"ahler manifolds with trivial canonical class}, Math. USSR-Sb.
{\bf 22} (1974) 580 - 583.

\bibitem[Bo2]{_Bogomolov:defo_} 
F. Bogomolov, {\em Hamiltonian K\"ahler
manifolds}, Sov. Math. Dokl. {\bf 19} (1978) 1462 - 1465.


%\bibitem[Bou]{_Boucksom_}
% Boucksom, S.,
%{\em  Higher dimensional Zariski decompositions},
% Ann. Sci. Ecole Norm. Sup. (4) 37 (2004), no. 1,
%   45--76, arXiv:math/0204336 


%\bibitem[Br]{_Brody:hyperboloc_}
% Brody, R., {\em Compact manifolds and hyperbolicity,}
%Trans. Amer. Math. Soc. 235 (1978), 213-219.

 \bibitem[BL1]{BL}
	Buzzard, G.; Lu, S. S.-Y.,
	{\em Algebraic surfaces holomorphically dominable by $\C^2$,}
	Invent. Math. 139, No. 3 (2000) 617-659. 

\bibitem[Camp]{_Camp04_}
        Campana, F.,
	{\em Orbifolds, special varieties and classification theory}, 
	{Ann. Inst. Fourier (Grenoble)}, {\bf 54} (2004), no.~3, 499 --630.


\bibitem[COP]{_COP:non-alge_}
Frederic Campana, Keiji Oguiso, Thomas Peternell,
{\em Non-algebraic hyperkaehler manifolds},
 J. Differential Geom. Volume 85, Number 3 (2010), 397-424. 

\bibitem[Cas]{_Cassels_}
        Cassels, J. W. S., 
        {\em Rational quadratic forms}, Academic Press, London, 1978. 

\bibitem[Cat]{_Catanese:moduli_}
F. Catanese, 
{\em A Superficial Working Guide to Deformations and Moduli},
arXiv:1106.1368, 
Advanced Lectures in Mathematics, Volume XXVI
Handbook of Moduli, Volume III, page 161-216 
(International Press).


\bibitem[D]{_Demailly:kobayashi_}
Demailly, Jean-Pierre,
{\em Hyperbolic algebraic varieties and holomorphic differential equations,}
 Talk given at the VIASM annual meeting, 25-26 August
 2012, long expository version (135 pages),
\url{http://www-fourier.ujf-grenoble.fr/~demailly/manuscripts/viasm2012expanded.pdf}.

\bibitem[FK]{FK} 
	H. M. Farkas and I. Kra, 
	{\em  Riemann Surfaces}, 
	Graduate Texts in Math. {71}, 
	Springer Verlag, New York, 1980.


\bibitem[F]{_Fujiki:HK_}  
Fujiki, A. {\em On the de Rham Cohomology Group of a Compact 
K\"ahler Symplectic Manifold}, Adv. Stud.
Pure Math. 10 (1987) 105-165.

\bibitem[GL]{GL} Greb, D.; Lehn, C., {\em Base manifolds for
Lagrangian fibrations on hyperk\"ahler manifolds}, 
arXiv:1303.3919v2


\bibitem[H1]{_Huybrechts:basic_} 
Huybrechts, D., 
{\em Compact hyperk\"ahler manifolds: Basic
results}, Invent. Math. 135 (1999), 63-113, 
alg-geom/9705025.


\bibitem[H2]{_Huybrechts:erratum_} 
Huybrechts, D., 
{\em Erratum to the paper: Compact hyperk\"ahler
  manifolds: basic results},
 Invent. math. 152  (2003), 209-212, math.AG/0106014.


%\bibitem[Hu3]{_Huybrechts:lec_}
% Huybrechts, Daniel, 
%{\em Compact hyperk\"ahler manifolds, Calabi-Yau 
%manifolds and related geometries,}
%Universitext, Springer-Verlag, Berlin, 2003, 
%Lectures from the Summer School held in
%Nordfjordeid, June 2001, pp. 161-225.

\bibitem[Hw0]{Hwang}
Jun-Muk Hwang, 
{\em Base manifolds for fibrations of projective 
irreducible symplectic manifolds,} Invent. Math. 174 (2008), no. 3, 
625--644.


\bibitem[KV]{_Kamenova_V:fibrations_}
Ljudmila Kamenova, Misha Verbitsky 
{\em Families of Lagrangian fibrations on hyperkaehler
  manifolds}, arXiv:1208.4626, 13 pages.

\bibitem[Kaw]{_Kawamata:pluricanonical_}
Y. Kawamata, {\em Pluricanonical systems on minimal algebraic varieties,} 
Invent. Math., {\bf 79}, (1985), no. 3, 567-588.


%\bibitem[Ko]{_Kobayashi:1970_}
%Kobayashi, Shoshichi,
%{\em Hyperbolic manifolds and holomorphic mappings,}
% Pure and Applied Mathematics 2,  (1970), New York: Marcel Dekker Inc.

\bibitem[Ko1]{_Kobayashi:1976_}
Kobayashi, Shoshichi,
{\em Intrinsic distances, measures and geometric function theory,}
Bulletin A.M.S. 82,  (1976), 357-416.

\bibitem[Ko2]{_Kobayashi:1998_}
S.~Kobayashi,
Hyperbolic complex spaces, Grundlehren der Mathematischen Wissenschaften,
Vol.~\textbf{318}, Springer-Verlag, Berlin, 1998.

%\bibitem[L]{_Lang:hyperbolic_}
%Lang, S., {\em 
%Introduction to complex hyperbolic spaces}, Springer, New York, 1987.

 \bibitem[LT]{LT}
Lang,  S.; Tate, J.,
 {\em Principal homogeneous spaces over abelian varieties},
 {American Journal of Mathematics.} 80, (1958), 659--684.

  \bibitem[Lu0]{_Lu:Multiply_}
Lu, S., {\em 
Multiply marked Riemann surfaces and the 
Kobayashi pseudometric
on algebraic manifolds,} Preprint 1999, arXiv:1308.1800

	\bibitem[Lu1]{_Lu:On_Kob_}
        Lu, S., {\em
        The Kobayashi pseudometric on algebraic manifold 
	and a canonical fibration,} Preprint 2001:
		arXiv:math.AG/0206170 and MPIM online preprint.
		
\bibitem[Lu2]{_Lu:Hol_} Lu, Steven S.-Y., 
{\em Holomorphic curves on irregular varieties of general type 
starting from surfaces}, Affine algebraic geometry, 205-220, CRM Proc. Lecture Notes, 54, Amer. Math. Soc., Providence, RI, 2011. 

\bibitem[Mar1]{_Markman:constra_}
Markman, E. {\em
Integral constraints on the monodromy group of
    the hyperkahler resolution of a symmetric product of a
    K3 surface,} International Journal of Mathematics
Vol. 21, No. 2 (2010) 169-223, arXiv:math/0601304.


%\bibitem[Mar2]{_Markman:survey_}
%Markman, E. {\em
%A survey of Torelli and monodromy results for 
%holomorphic-symplectic varieties}, 
%Proceedings of the conference "Complex and Differential 
%Geometry'', Springer Proceedings in Mathematics, 2011, Volume 8, 257--322,
%arXiv:math/0601304.

%\bibitem[Mar3]{_Markman:divisors_}
%Markman, E. {\em
%Prime exceptional divisors on holomorphic symplectic varieties and 
%monodromy-reflections}, 
%arXiv: 0912.4981, 53 pages

\bibitem[Mar4]{_Markman:SYZ_}
Markman, E. 
{\em Lagrangian fibrations of holomorphic-symplectic varieties of 
$K3^{[n]}$-type}, arXiv:1301.6584 [math.AG]

\bibitem[Mat1]{_Matsushita:fibred_} 
D. Matsushita, {\em On fibre space structures of a
  projective irreducible symplectic manifold},
alg-geom/9709033, math.AG/9903045, also in Topology
{\bf 38} (1999), No. 1, 79-83. Addendum, 
Topology {\bf 40} (2001) No. 2, 431 - 432.

%\bibitem[Mat2]{_Matsushita:CP^n_} 
% Matsushita, D., {\em
%Higher direct images of Lagrangian fibrations},
%Amer. J. Math. 127 (2005) arXiv:math/0010283.

\bibitem[MZ]{_Matsushita_Zhang_} 
 Matsushita, D.; Zhang, De-Qi, {\em
Zariski F-decomposition and Lagrangian fibration on 
hyperk\"ahler manifolds}, arXiv:0907.5311v4.

\bibitem[Me]{_Meyers_}
A. Meyer, {\em Mathematische Mittheilungen}, 
Vierteljahrschrift der Naturforschenden Gesellschaft in Z\"urich 29 (1884), 209-22.

\bibitem[MO]{_Mongardi_Onorati_} 
Mongardi, G., Onorati, C., 
{\em Birational geometry of irreducible holomorphic symplectic tenfolds 
of O'Grady type}, 
arXiv:2010.12511 [math.AG]. 

\bibitem[MR]{_Mongardi_Rapagnetta_}
Mongardi, G., Rapagnetta, A., 
{\em Monodromy and birational geometry of O'Grady's sixfolds}, 
arXiv:1909.07173 [math.AG]. 

 \bibitem[MM]{MM} Mori, S.;  Mukai, S.,
{\em The uniruledness of the moduli space of curves of genus 11,}
	Lecture Notes in Math.
	1016 (1982), 334--353.
	
 \bibitem[Roy]{Roy}
	Royden, H.,
	{\em The extension of regular holomorphic maps,}
	Proc. Amer. Math. Soc. 43 (1974) 306--310.

%\bibitem[PS]{_Piatetski_Shapiro_Shafarevich_}
%Piatecki-Shapiro, I.I.; Shafarevich I.R.,
%{\em Torelli's
%theorem for algebraic surfaces of type
%K3}, Izv. Akad. Nauk SSSR Ser. Mat. (1971) 35: 530-572.


\bibitem[Saw]{_Sawon_} 
Sawon, J.,  
{\em Abelian fibred holomorphic symplectic
  manifolds}, Turkish
Jour. Math. 27 (2003) no. 1, 197 - 230, math.AG/0404362.

\bibitem[Siu]{Siu} Siu, Y.-T., {\em Every Stein subvariety admits
a Stein neighbourhood},
Inventiones Math. 38 (1976), 89 - 100.

%\bibitem[Su]{_Sullivan:infinite_}
%Sullivan, D.,
%{\em Infinitesimal computations in topology}, Publications
%Math\-\'ema\-tiques de l'IH\'ES, 47 (1977), p. 269-331


\bibitem[V0]{_Verbitsky:coho_announce_} 
Verbitsky, M.,
{\it Cohomology of compact hyperk\"ahler manifolds
and its applications,}  GAFA vol. 6 (4) pp. 601-612 (1996).


\bibitem[V1]{_Verbitsky:SYZ_}
Verbitsky, M., 
{\em Hyperk\"ahler SYZ conjecture and semipositive line bundles},
 arXiv:0811.0639, 21 pages, GAFA 19, No. 5 (2010) 1481-1493.


\bibitem[V2]{_V:Torelli_}
Verbitsky, M.,
{\em A global Torelli theorem for hyperk\"ahler manifolds,}
arXiv: 0908.4121, 52 pages (to appear in Duke Math. J.)

\bibitem[V3]{_Verbitsky:ergodic_}
Verbitsky, M.,
{\em Ergodic complex structures on hyperk\"ahler manifolds},
arXiv:1306.1498, 22 pages.


\bibitem[Vo]{_Voisin:kobayashi_}
Claire Voisin,
{\em On some problems of Kobayashi and Lang; algebraic
approaches,}  Current Developments in Mathematics 2003,
no. 1 (2003), 53-125.


%\bibitem[Vi]{_Viehweg:moduli_}
%Viehweg, E.,
%{\em Quasi-projective Moduli for Polarized Manifolds,}
%Springer-Verlag, Berlin, Heidelberg, New York, 1995,
%Ergebnisse der Mathematik und ihrer Grenzgebiete, 3.~Folge, Band 30,
%also available at
%{\tt http://www.uni-due.de/$\widetilde{\phantom{a}}$mat903/books.html}


%\bibitem[VGS]{_Vinberg_Gorbatsevich_Shvartsman_}
%Vinberg, E. B.,  Gorbatsevich, V. V.,  Shvartsman, O. V., 
%{\em Discrete Subgroups of Lie Groups}, in 
%``Lie Groups and Lie Algebras II'', Springer-Verlag, 2000.

 \bibitem[Yau]{Yau}
	Yau, S.T.,
	{\em On the Ricci-curvature of a complex K\"ahler manifold and
	the complex Monge-Amp\`ere equation},
	{Comm. Pure Appl. Math.},
	31 (1978), 339--411.

\bibitem[Y]{_Yoshioka_} Yoshioka, K., 
        {\em Bridgeland's stability and the positive cone of the moduli spaces 
        of stable objects on an abelian surface}, 
        arXiv:1206.4838 [math.AG]. 

\bibitem[Zai]{Zai} Zaidenberg, M. G., {\em
On hyperbolic embeddings of complements of divisors and the limiting behavior of the Kobayashi-Royden metric}, 
Math. USSR-Sb. 55 (1986), 55 - 70.
\end{thebibliography}
\end{document}